# Determination of the Stieltjes constants at rational arguments

Donal F. Connon

19 June 2015

## Abstract

The Stieltjes constants $\gamma_n(x)$ have attracted considerable attention in recent years and a number of authors, including the present one, have considered various ways in which these constants may be evaluated. The primary purpose of this paper is to belatedly highlight the fact that Deninger [15] actually ascertained the first generalised Stieltjes constant at rational arguments as long ago as 1984 and that all of the higher constants (at rational arguments) were determined in principle by Chakraborty, Kanemitsu and Kuzumaki [8] in 2009. Equivalent results were obtained by Musser [21] in his 2011 thesis.

The authors of the papers [15] and [8] simply referred to the $\gamma_n(x)$ as the Laurent coefficients which explains why various electronic searches conducted by this author for "Stieltjes constants" did not readily highlight these particular sources.

In this paper the author has employed a slightly different argument to obtain a simpler expression for the results originally derived by Chakraborty et al. [8] in 2009.

## 1. Introduction

The Stieltjes constants $\gamma_n(x)$ are the coefficients of the Laurent expansion of the Hurwitz zeta function $\varsigma(s,x)$ about $s=1$

$$(1.1) \qquad \varsigma(s,x)=\sum_{n=0}^{\infty}\frac{1}{(n+x)^s}=\frac{1}{s-1}+\sum_{n=0}^{\infty}\frac{(-1)^n}{n!}\gamma_n(x)(s-1)^n$$

where $\gamma_p(x)$ are known as the generalised Stieltjes constants and we have [25]

$$(1.2) \qquad \gamma_0(x)=-\psi(x)$$

where $\psi(x)$ is the digamma function.

With $x=1$ equation (1.1) reduces to the Laurent expansion of the Riemann zeta function

$$\varsigma(s)=\frac{1}{s-1}+\sum_{n=0}^{\infty}\frac{(-1)^n}{n!}\gamma_n(s-1)^n$$

where $\gamma_n\equiv\gamma_n(1)$.

## 2. A useful formula for the Stieltjes constants

We recall Hasse's formula [16] for the Hurwitz zeta function which is valid for all $s \in \mathbf{C}$ except $s = 1$ (it is valid in the limit as $s \to 1$ in the form shown below)

(2.1) $$(s-1)\varsigma(s,x) = \sum_{j=0}^{\infty} \frac{1}{j+1} \sum_{k=0}^{j} \binom{j}{k} \frac{(-1)^k}{(x+k)^{s-1}}$$

and differentiation with respect to $x$ gives us

$$\frac{\partial}{\partial x}\varsigma(s,x) = -\sum_{j=0}^{\infty} \frac{1}{j+1} \sum_{k=0}^{j} \binom{j}{k} \frac{(-1)^k}{(x+k)^{s}}$$

and we then have

(2.2) $$\frac{\partial^n}{\partial s^n}\frac{\partial}{\partial x}\varsigma(s,x) = (-1)^{n+1} \sum_{j=0}^{\infty} \frac{1}{j+1} \sum_{k=0}^{j} \binom{j}{k} \frac{(-1)^k \log^n(x+k)}{(x+k)^{s}}$$

We note that the partial derivatives commute in the region where $\varsigma(s,x)$ is analytic and hence we have

$$\frac{\partial^n}{\partial s^n}\frac{\partial}{\partial x}\varsigma(s,x) = \frac{\partial}{\partial x}\frac{\partial^n}{\partial s^n}\varsigma(s,x)$$

Evaluation of (2.2) at $s = 0$ results in

(2.3) $$\frac{\partial}{\partial x}\varsigma^{(n)}(0,x) = (-1)^{n+1} \sum_{j=0}^{\infty} \frac{1}{j+1} \sum_{k=0}^{j} \binom{j}{k} (-1)^k \log^n(x+k)$$

The following expression for the Stieltjes constants was previously obtained by the author in [12] in 2007

(2.4) $$\gamma_n(x) = -\frac{1}{n+1} \sum_{j=0}^{\infty} \frac{1}{j+1} \sum_{k=0}^{j} \binom{j}{k} (-1)^k \log^{n+1}(x+k)$$

and hence, as was also shown by Chakraborty, Kanemitsu and Kuzumaki [8], we have deduced the important identity

(2.5) $$R_n'(x) = (-1)^{n+1} \frac{\partial^n}{\partial s^n}\frac{\partial}{\partial x}\varsigma(s,x)\bigg|_{s=0} = -n\gamma_{n-1}(x)$$

where for $n \geq 1$

$$R_n(x) = (-1)^{n+1} \frac{\partial^n}{\partial s^n}\varsigma(s,x)\bigg|_{s=0}$$

In Eq. (4.3.231) of [12] we also noted that $R_n'(x) = -n\gamma_{n-1}(x)$ in the equivalent form for $x > 0$

$$\int_1^x \gamma_n(t)\,dt = \frac{(-1)^{n+1}}{n+1}\left[\varsigma^{(n+1)}(0,x) - \varsigma^{(n+1)}(0)\right]$$

but the usefulness of this simple formula was not then fully appreciated by this author.

Letting $s \to 1-s$ in (1.1) gives us

$$\varsigma(1-s,x) = -\frac{1}{s} + \sum_{n=0}^{\infty} \frac{1}{n!}\gamma_n(x)s^n$$

and substituting $R_{n+1}'(x) = -(n+1)\gamma_n(x)$ we obtain

$$\varsigma(1-s,x) = -\frac{1}{s} - \sum_{n=0}^{\infty} \frac{R_{n+1}'(x)}{(n+1)!}s^n$$

Reindexing gives us

$$\varsigma(1-s,x) = -\frac{1}{s} - \frac{1}{s}\sum_{n=1}^{\infty} \frac{R_n'(x)}{n!}s^n$$

and thus we have

(2.6) $$s\varsigma(1-s,x) = -1 - \sum_{n=1}^{\infty} \frac{R_n'(x)}{n!}s^n$$

The Deninger $R-$function is defined by

$$R(x) = R_2(x) = -\varsigma''(0,x)$$

and it is named after Deninger [15] who introduced it in 1984.

## 3. Determination of the Stieltjes constants at rational arguments

We recall Rademacher's formula [1, p.261] for the Hurwitz zeta function where for all $s$ and $1 \le p \le q$ where $p$ and $q$ are positive integers

(3.1) $$\varsigma\left(s, \frac{p}{q}\right) = 2\Gamma(1-s)(2\pi q)^{s-1}\sum_{j=1}^{q} \sin\left(\frac{\pi s}{2} + \frac{2\pi jp}{q}\right)\varsigma\left(1-s, \frac{j}{q}\right)$$

and letting $s \to 1-s$ gives us

(3.2) $$\varsigma\left(1-s, \frac{p}{q}\right) = 2\Gamma(s)(2\pi q)^{-s}\sum_{j=1}^{q} \cos\left(\frac{2\pi jp}{q} - \frac{\pi s}{2}\right)\varsigma\left(s, \frac{j}{q}\right)$$

Equating this with (2.6) gives us

(3.3) $$2\Gamma(1+s)(2\pi q)^{-s}\sum_{r=1}^{q}\cos\left(\frac{2\pi rp}{q}-\frac{\pi s}{2}\right)\varsigma\left(s,\frac{r}{q}\right)=-1-\sum_{n=1}^{\infty}R_n'\left(\frac{p}{q}\right)\frac{s^n}{n!}$$

For convenience we designate $\phi_r(s)$ as

$$\phi_r(s)=\Gamma(1+s)(2\pi q)^{-s}\cos\left(\frac{2\pi rp}{q}-\frac{\pi s}{2}\right)\varsigma\left(s,\frac{r}{q}\right)\equiv f_1(s)f_2(s)f_3(s)f_4(s)$$

and using the Leibniz rule for derivatives we have with $D^n\equiv\frac{d^n}{ds^n}$

$$D^n\phi_r(s)=\sum_{i=0}^{n}\binom{n}{i}D^i[f_1(s)f_2(s)]D^{n-i}[f_3(s)f_4(s)]$$

where

$$D^{n-i}[f_3(s)f_4(s)]=\sum_{j=0}^{n-i}\binom{n-i}{j}D^{n-i-j}f_3(s)D^jf_4(s)$$

We see that

(3.4) $$\frac{d}{ds}\Gamma(1+s)(2\pi q)^{-s}=\Gamma(1+s)(2\pi q)^{-s}[\psi(1+s)-\log(2\pi q)]$$

and hence, as explained below, we may present the higher derivatives in terms of the (exponential) complete Bell polynomials which may be defined by $Y_0=1$ and for $r\geq 1$

$$Y_r(x_1,...,x_r)=\sum_{\pi(r)}\frac{r!}{k_1!\,k_2!...\,k_r!}\left(\frac{x_1}{1!}\right)^{k_1}\left(\frac{x_2}{2!}\right)^{k_2}...\left(\frac{x_r}{r!}\right)^{k_r}$$

where the sum is taken over all partitions $\pi(r)$ of $r$, i.e. over all sets of integers $k_j$ such that

$$k_1+2k_2+3k_3+\cdots+rk_r=r$$

The complete Bell polynomials have integer coefficients and the first six are set out below (Comtet [11, p.307])

(3.5) $$Y_1(x_1)=x_1$$

$$Y_2(x_1,x_2)=x_1^2+x_2$$

$$Y_3(x_1,x_2,x_3)=x_1^3+3x_1x_2+x_3$$

$$Y_4(x_1,x_2,x_3,x_4)=x_1^4+6x_1^2x_2+4x_1x_3+3x_2^2+x_4$$

$$Y_5(x_1,x_2,x_3,x_4,x_5) = x_1^5 + 10x_1^3 x_2 + 10x_1^2 x_3 + 15x_1 x_2^2 + 5x_1 x_4 + 10x_2 x_3 + x_5$$

$$Y_6(x_1,x_2,x_3,x_4,x_5,x_6) = x_1^6 + 6x_1 x_5 + 15x_2 x_4 + 10x_2^3 + 15x_1^2 x_4 + 15x_2^3 + 60x_1 x_2 x_3$$

$$+20x_1^3 x_3 + 45x_1^2 x_2^2 + 15x_1^4 x_1 + x_6$$

For our purposes, all of the information that we require regarding these polynomials is contained in [14] and more detailed expositions may be found in [9] and [11].

In particular, suppose that $h'(x) = h(x)g(x)$ then we have

(3.6) $$h^{(r)}(x) = h(x) Y_r\left(g(x), g^{(1)}(x), ..., g^{(r-1)}(x)\right)$$

Hence, having regard to (3.4) and (3.6), we obtain

$$D^i \Gamma(1+s)(2\pi q)^{-s} = \Gamma(1+s)(2\pi q)^{-s} Y_i\left(\psi(1+s) - \log(2\pi q), \psi'(1+s), \cdots, \psi^{(i-1)}(1+s)\right)$$

and, for convenience, we designate $\mathbf{Y}_i\left(\boldsymbol{\psi}(s)\right)$ as

$$\mathbf{Y}_i\left(\boldsymbol{\psi}(s)\right) = Y_i\left(\psi(1+s) - \log(2\pi q), \psi'(1+s), \cdots, \psi^{(i-1)}(1+s)\right)$$

We then have

$$D^i \Gamma(1+s)(2\pi q)^{-s}\Big|_{s=0} = \mathbf{Y}_i\left(\boldsymbol{\psi}(0)\right)$$

where $\mathbf{Y}_i\left(\boldsymbol{\psi}(0)\right) = Y_i\left(-[\gamma + \log(2\pi q)], \psi'(1), \cdots, \psi^{(i-1)}(1)\right)$.

It is well known that [24, p.22]

$$\psi^{(p)}(x) = (-1)^{p+1} p!\varsigma(p+1, x)$$

so that

(3.7) $$\psi^{(p)}(1) = (-1)^{p+1} p!\varsigma(p+1,1) = (-1)^{p+1} p!\varsigma(p+1)$$

We then have

$$\mathbf{Y}_i\left(\boldsymbol{\psi}(0)\right) = Y_i\left(-[\gamma + \log(2\pi q)], 1!\varsigma(2), \cdots, (-1)^i (i-1)!\varsigma(i)\right)$$

From the definition of the (exponential) complete Bell polynomials we have

$$Y_m(ax_1, a^2 x_2, ..., a^m x_m) = a^m Y_r(x_1, ..., x_m)$$

and thus with $a = -1$ we have

(3.8) $$Y_m(-x_1, x_2, ..., (-1)^m x_m) = (-1)^m Y_m(x_1, ..., x_m)$$

Since $Y_m(x_1, ..., x_m) > 0$ when all of the arguments are positive, we deduce that $Y_m(-x_1, x_2, ..., (-1)^m x_m)$ has the same sign as $(-1)^m$.

Accordingly, we may write

$$\mathbf{Y}_i\left(\boldsymbol{\psi}(0)\right) = (-1)^i Y_i\left(\gamma + \log(2\pi q), 1!\varsigma(2), \cdots, (i-1)!\varsigma(i)\right)$$

which we abbreviate to $\mathbf{Y}_i\left(\boldsymbol{\psi}(0)\right) = (-1)^i \mathbf{Y}_i^*\left(\boldsymbol{\psi}(0)\right)$.

With $f_3(s) = \cos\left(\frac{2\pi rp}{q} - \frac{\pi s}{2}\right)$ we see that

$$f_3'(s) = \frac{\pi}{2}\sin\left(\frac{2\pi rp}{q} - \frac{\pi s}{2}\right) = -\frac{\pi}{2}\cos\left(\frac{2\pi rp}{q} - \frac{\pi s}{2} + \frac{\pi}{2}\right)$$

and hence we obtain

$$D^{n-i-j} f_3(s) = (-1)^{n-i-j}\left(\frac{\pi}{2}\right)^{n-i-j}\cos\left(\frac{2\pi rp}{q} - \frac{\pi s}{2} + \frac{\pi(n-i-j)}{2}\right)$$

We therefore have

$$D^{n-i-j} f_3(s)\Big|_{s=0} = (-1)^{n-i-j}\left(\frac{\pi}{2}\right)^{n-i-j}\cos\left(\frac{2\pi rp}{q} + \frac{\pi(n-i-j)}{2}\right)$$

and by definition we have

$$D^j f_4(s)\Big|_{s=0} = \varsigma^{(j)}\left(0, \frac{r}{q}\right)$$

This gives us

$$D^n \phi_r(s)\Big|_{s=0} = \sum_{i=0}^{n}\binom{n}{i}\mathbf{Y}_i^*\left(\boldsymbol{\psi}(0)\right)\sum_{j=0}^{n-i}\binom{n-i}{j}(-1)^{n-j}\left(\frac{\pi}{2}\right)^{n-i-j}\cos\left(\frac{2\pi rp}{q} + \frac{\pi(n-i-j)}{2}\right)\varsigma^{(j)}\left(0, \frac{r}{q}\right)$$

From (3.3) we see that

(3.9) $$2\sum_{r=1}^{q} D^n \phi_r(s)\Big|_{s=0} = -R_n'\left(\frac{p}{q}\right) = n\gamma_{n-1}\left(\frac{p}{q}\right)$$

and hence we obtain

(3.10)

$$n\gamma_{n-1}\left(\frac{p}{q}\right)=2\sum_{r=1}^{q}\sum_{i=0}^{n}\binom{n}{i}\mathbf{Y}_i^*\left(\boldsymbol{\psi}(0)\right)\sum_{j=0}^{n-i}\binom{n-i}{j}(-1)^{n-j}\left(\frac{\pi}{2}\right)^{n-i-j}\cos\left(\frac{2\pi rp}{q}+\frac{\pi(n-i-j)}{2}\right)\varsigma^{(j)}\left(0,\frac{r}{q}\right)$$

A similar approach may be found in Musser's thesis [21, p.25].

For example, in the simplest case with $n=1$, after some tedious algebra we obtain

$$\gamma_0\left(\frac{p}{q}\right)=\gamma+\log(2\pi q)+2\sum_{r=1}^{q}\left[\cos\left(\frac{2\pi rp}{q}\right)\varsigma'\left(0,\frac{r}{q}\right)+\frac{\pi}{2}\sin\left(\frac{2\pi rp}{q}\right)\varsigma\left(0,\frac{r}{q}\right)\right]$$

We have the well-known identity [1, p.264]

(3.11) $$\varsigma(0,x)=\frac{1}{2}-x$$

Using Lerch's identity [24, p.92] for $x>0$

(3.12) $$\varsigma'(0,x)=\log\Gamma(x)-\frac{1}{2}\log(2\pi)$$

and the trigonometric formulae

$$\sum_{r=1}^{q}\sin\left(\frac{2\pi rp}{q}\right)=0$$

$$\sum_{r=1}^{q}\cos\left(\frac{2\pi rp}{q}\right)=0 \qquad p=1,2,...,q-1$$

$$\sum_{r=1}^{q}r\sin\left(\frac{2\pi rp}{q}\right)=-\frac{q}{2}\cot\left(\frac{p\pi}{q}\right) \quad \text{for } p<q$$

we obtain using (1.2)

(3.14) $$\psi\left(\frac{p}{q}\right)=-\gamma-\log(2\pi q)-\frac{\pi}{2}\cot\left(\frac{p\pi}{q}\right)-2\sum_{r=1}^{q}\cos\left(\frac{2\pi rp}{q}\right)\log\Gamma\left(\frac{r}{q}\right)$$

which concurs with the format for the Gauss identity for the digamma function noted by Ram Murty and Saradhab [22].

Euler's reflection formula for the gamma function gives us

$$-\sum_{r=1}^{q-1}\cos\left(\frac{2\pi rp}{q}\right)\log\sin\left(\frac{\pi r}{q}\right)=\sum_{r=1}^{q-1}\cos\left(\frac{2\pi rp}{q}\right)\left[\log\Gamma\left(\frac{r}{q}\right)+\log\Gamma\left(1-\frac{r}{q}\right)-\log\pi\right]$$

$$= \log \pi + \sum_{r=1}^{q-1} \cos\left(\frac{2\pi rp}{q}\right)\left[\log\Gamma\left(\frac{r}{q}\right) + \log\Gamma\left(1-\frac{r}{q}\right)\right]$$

and we note that

$$\sum_{r=1}^{q-1} \cos\left(\frac{2\pi rp}{q}\right)\log\Gamma\left(1-\frac{r}{q}\right) = \sum_{r=1}^{q-1} \cos\left[2\pi p\left(1-\frac{r}{q}\right)\right]\log\Gamma\left(1-\frac{r}{q}\right)$$

$$= \sum_{r=1}^{q-1} \cos\left[2\pi p\left(\frac{q-r}{q}\right)\right]\log\Gamma\left(\frac{q-r}{q}\right)$$

With new indexing, $r = q - j$, the latter summation becomes $\sum_{j=1}^{q-1} \cos\left(\frac{2\pi jp}{q}\right)\log\Gamma\left(\frac{j}{q}\right)$ and hence we deduce that

$$-\sum_{r=1}^{q-1} \cos\left(\frac{2\pi rp}{q}\right)\log\sin\left(\frac{\pi r}{q}\right) = \log\pi + 2\sum_{r=1}^{q-1} \cos\left(\frac{2\pi rp}{q}\right)\log\Gamma\left(\frac{r}{q}\right)$$

Thus we may write (3.14) as [15] for $p < q$

(3.15) $$\psi\left(\frac{p}{q}\right) = -\gamma - \log(2q) - \frac{\pi}{2}\cot\left(\frac{p\pi}{q}\right) + \sum_{r=1}^{q-1} \cos\left(\frac{2\pi rp}{q}\right)\log\sin\left(\frac{\pi r}{q}\right)$$

I drafted a detailed paper on the Stieltjes constants in 2011 and it was then that I came across [8] and rather belatedly ascertained that, as long ago as 1984, Deninger [15] had determined that (see also Eq. (2.2) of [8])

(3.16) $$R'\left(\frac{p}{q}\right) = -2\gamma_1 - 2[\gamma + \log(2\pi q)]\left[\gamma + \psi\left(\frac{p}{q}\right)\right] + 2\sum_{j=1}^{q-1} \cos\left(\frac{2\pi jp}{q}\right)R\left(\frac{j}{q}\right)$$

$$-\pi\sum_{j=1}^{q-1} \sin\left(\frac{2\pi jp}{q}\right)\log\Gamma\left(\frac{j}{q}\right) - \log^2 q - 2\log q\log(2\pi)$$

where $R(x) = -\varsigma''(0,x)$ and $R'(x) = -2\gamma_1(x)$. In the above, I have corrected a typographical error in the term involving $\sum_{j=1}^{q-1} \sin\left(\frac{2\pi jp}{q}\right)\log\Gamma\left(\frac{j}{q}\right)$.

The discovery of this earlier work made most of my 2011 paper no longer germane and hence I did not publish it (Blagouchine [7] kindly referred to some of the results from this draft paper in his extensive 2015 paper relating to the Stieltjes constants).

It should be noted that the Stieltjes constant $\gamma_1$ reported in Eq. (2.15) of [15] should be increased by a factor of 2 (this difference arises because Deninger [15, p.174] has employed a different definition in the Laurent expansion of the Hurwitz zeta function; this value was also inconsistently employed in Eq. (1.26) of the paper by Chakraborty, Kanemitsu, and Kuzumaki [8]).

We may write (3.16) as

$$(3.17) \qquad \gamma_1\left(\frac{p}{q}\right)=\gamma_1+[\gamma+\log(2\pi q)]\left[\gamma+\psi\left(\frac{p}{q}\right)\right]+\sum_{j=1}^{q-1}\cos\left(\frac{2\pi jp}{q}\right)\varsigma''\left(0,\frac{j}{q}\right)$$

$$+\pi\sum_{j=1}^{q-1}\sin\left(\frac{2\pi jp}{q}\right)\log\Gamma\left(\frac{j}{q}\right)+\frac{1}{2}\log^2 q+\log q\log(2\pi)$$

The formula for $R_n'\left(\frac{p}{q}\right)$ was generalised in 2009 by Chakraborty et al. [8] to give a corresponding expression for $\gamma_n\left(\frac{p}{q}\right)$. Their representation is recorded below:

$$(3.18) \qquad R_n'\left(\frac{p}{q}\right)=2(-1)^n\sum_{l=0}^{n}\binom{n}{l}\delta_{n-l}\sum_{r=1}^{q}K_l\left(\frac{p}{q},r\right)$$

where:

$$-K_l\left(\frac{p}{q},r\right)=\sum_{m=0}^{l}\binom{l}{m}R_{l-m}\left(\frac{r}{q}\right)\left(\frac{\pi}{2}\right)^m T_m\left(\frac{rp}{q}\right)$$

$$T_m=\begin{cases}(-1)^{m/2}\cos(2\pi x), & m \text{ even}\\ (-1)^{(m+1)/2}\sin(2\pi x), & m \text{ odd}\end{cases}$$

$$\delta_l=\sum_{m=0}^{l}(-1)^m\binom{l}{m}\Gamma^{(m)}(1)\log^{l-m}(2\pi q)$$

Prima facie, the formula above is more complex than (3.10) in that it contains the derivatives of the gamma function $\Gamma^{(m)}(1)$ whereas the representation (3.10) in this paper does not require that additional level of complexity; instead we only need to consider the simpler expression $\mathbf{Y}_i^*\left(\boldsymbol{\psi}(0)\right)=Y_i\left(\gamma+\log(2\pi q),1!\varsigma(2),\cdots,(i-1)!\varsigma(i)\right)$ which basically just employs the values $\gamma$, $\log(2\pi q)$ and $\varsigma(i)$.

Formula (3.18) contains four nested summations whereas there are only two summations in (3.10), albeit the latter expression contains Bell polynomials. Further simplifications to the latter formula are considered below:

We have

$$(3.19) \qquad \phi_r(s)=\Gamma(1+s)(2\pi q)^{-s}\cos\left(\frac{\pi s}{2}\right)\varsigma\left(s,\frac{r}{q}\right)\cos\left(\frac{2\pi rp}{q}\right)$$

$$+\Gamma(1+s)(2\pi q)^{-s}\sin\left(\frac{\pi s}{2}\right)\varsigma\left(s,\frac{r}{q}\right)\sin\left(\frac{2\pi rp}{q}\right)$$

and using Euler's reflection formula

$$\Gamma(s)\Gamma(1-s)=\frac{\pi}{\sin\pi s}$$

we may express (3.19) as

(3.20) $$\phi_r(s)=(2\pi q)^{-s}\frac{\Gamma(1+s/2)\Gamma(1-s/2)}{\Gamma(1-s)}\varsigma\left(s,\frac{r}{q}\right)\cos\left(\frac{2\pi rp}{q}\right)$$

$$+\pi(2\pi q)^{-s}\frac{\Gamma(1+s)}{\Gamma(s/2)\Gamma(1-s/2)}\varsigma\left(s,\frac{r}{q}\right)\sin\left(\frac{2\pi rp}{q}\right)$$

$$\equiv\lambda(s)\varsigma\left(s,\frac{r}{q}\right)\cos\left(\frac{2\pi rp}{q}\right)+\mu(s)\varsigma\left(s,\frac{r}{q}\right)\sin\left(\frac{2\pi rp}{q}\right)$$

where we have defined $\lambda(s)$ as

$$\lambda(s)=(2\pi q)^{-s}\frac{\Gamma(1+s/2)\Gamma(1-s/2)}{\Gamma(1-s)}$$

We see that

$$\lambda'(s)=\lambda(s)g(s)$$

where

$$g(s)=\psi(1-s)+\frac{1}{2}\psi\left(1+\frac{s}{2}\right)-\frac{1}{2}\psi\left(1-\frac{s}{2}\right)-\log(2\pi q)$$

Therefore, in accordance with (3.6), the higher derivatives may be expressed as

$$\lambda^{(i)}(s)=\lambda(s)Y_i\left(g(s),g^{(1)}(s),...,g^{(i-1)}(s)\right)$$

where

$$g^{(i)}(s)=(-1)^i\psi^{(i)}(1-s)+\frac{1}{2^{i+1}}\psi^{(i)}\left(1+\frac{s}{2}\right)-\frac{1}{2^{i+1}}(-1)^i\psi^{(i)}\left(1-\frac{s}{2}\right)$$

We have

$$g^{(i)}(0)=\left[(-1)^i+\frac{1}{2^{i+1}}[1-(-1)^i]\right]\psi^{(i)}(1)$$

and substituting

$$\psi^{(i)}(1)=(-1)^{i+1}i!\varsigma(i+1,1)=(-1)^{i+1}i!\varsigma(i+1)$$

we obtain

$$g^{(i)}(0) = \left[ \frac{1}{2^{i+1}}[(-1)^{i+1}+1]-1 \right] i!\varsigma(i+1)$$

We have

$$\lambda^{(i)}(0) = Y_i\left(g(0), g^{(1)}(0), ..., g^{(i-1)}(0)\right) \equiv \mathbf{Y}_i(\mathbf{g}(0))$$

where, for example, we have

$$g(0) = -[\gamma + \log(2\pi q)]$$

$$g^{(1)}(0) = -\frac{1}{2}\varsigma(2)$$

We now consider the second limb of (3.20) and, for convenience, we have defined $\mu(s)$ as

$$\mu(s) = \pi(2\pi q)^{-s} \frac{\Gamma(1+s)}{\Gamma(s/2)\Gamma(1-s/2)}$$

$$= \frac{\pi}{2} s\,(2\pi q)^{-s} \frac{\Gamma(1+s)}{\Gamma(1+s/2)\Gamma(1-s/2)} \equiv \frac{\pi}{2} s\,\theta(s)$$

Hence we have

$$\mu^{(i)}(s) = \frac{\pi}{2}\left[ s\theta^{(i)}(s) + i\theta^{(i-1)}(s) \right]$$

and thus

$$\mu^{(i)}(0) = \frac{\pi}{2} i\theta^{(i-1)}(0)$$

It is easily seen that

$$\theta'(s) = \theta(s)h(s)$$

where

$$h(s) = \psi(1+s) - \frac{1}{2}\psi\left(1+\frac{s}{2}\right) + \frac{1}{2}\psi\left(1-\frac{s}{2}\right) - \log(2\pi q)$$

and hence we have

$$\theta^{(i-1)}(s) = \theta(s)Y_{i-1}\left(h(s), h^{(1)}(s), ..., h^{(i-2)}(s)\right)$$

where

$$h^{(i)}(s) = \psi^{(i)}(1+s) - \frac{1}{2^{i+1}}\psi^{(i)}\left(1+\frac{s}{2}\right) + \frac{1}{2^{i+1}}(-1)^i\psi^{(i)}\left(1-\frac{s}{2}\right)$$

It may be noted that $h(s) = g(-s)$, but this does not appear to simplify the analysis. We have

$$\mu^{(i)}(0) = \frac{\pi}{2} i Y_{i-1}\left(h(0), h^{(1)}(0), ..., h^{(i-2)}(0)\right)$$

and recalling (3.9)

$$2\sum_{r=1}^{q} D^n \phi_r(s)\Big|_{s=0} = -R_n'\left(\frac{p}{q}\right) = n\gamma_{n-1}\left(\frac{p}{q}\right)$$

we see that

(3.21) $$n\gamma_{n-1}\left(\frac{p}{q}\right) = 2\sum_{r=1}^{q}\sum_{i=0}^{n}\binom{n}{i}\left[\mathbf{Y}_i(\mathbf{g}(0))\cos\left(\frac{2\pi rp}{q}\right) + \frac{\pi}{2} i\mathbf{Y}_{i-1}(\mathbf{h}(0))\sin\left(\frac{2\pi rp}{q}\right)\right]\varsigma^{(n-i)}\left(0, \frac{r}{q}\right)$$

where $\mathbf{Y}_i(\mathbf{g}(0)) \equiv Y_i\left(g(0), g^{(1)}(0), ..., g^{(i-1)}(0)\right)$ and $\mathbf{Y}_{i-1}(\mathbf{h}(0)) \equiv Y_{i-1}\left(h(0), h^{(1)}(0), ..., h^{(i-2)}(0)\right)$.

With $p = q = 1$ we obtain

(3.22) $$n\gamma_{n-1} = 2\sum_{i=0}^{n}\binom{n}{i}\mathbf{Y}_i(\mathbf{g}(0))\varsigma^{(n-i)}(0)$$

which is also derived in a different manner in (5.4) below.

We may write (3.21) as

$$n\gamma_{n-1}\left(\frac{p}{q}\right) = 2\sum_{r=1}^{q-1}\sum_{i=0}^{n}\binom{n}{i}\left[\mathbf{Y}_i(\mathbf{g}(0))\cos\left(\frac{2\pi rp}{q}\right) + i\frac{\pi}{2}\mathbf{Y}_{i-1}(\mathbf{h}(0))\sin\left(\frac{2\pi rp}{q}\right)\right]\varsigma^{(n-i)}\left(0, \frac{r}{q}\right)$$

$$+2\sum_{i=0}^{n}\binom{n}{i}\mathbf{Y}_i(\mathbf{g}(0))\varsigma^{(n-i)}(0)$$

and hence we have

(3.23)
$$n\gamma_{n-1}\left(\frac{p}{q}\right) = n\gamma_{n-1} + 2\sum_{r=1}^{q-1}\sum_{i=0}^{n}\binom{n}{i}\left[\mathbf{Y}_i(\mathbf{g}(0))\cos\left(\frac{2\pi rp}{q}\right) + \frac{\pi}{2} i\mathbf{Y}_{i-1}(\mathbf{h}(0))\sin\left(\frac{2\pi rp}{q}\right)\right]\varsigma^{(n-i)}\left(0, \frac{r}{q}\right)$$

With $p = 1$ and $q = 2$ we have

(3.24) $$n\gamma_{n-1}\left(\frac{1}{2}\right) = n\gamma_{n-1} - 2\sum_{i=0}^{n}\binom{n}{i}\mathbf{Y}_i(\mathbf{g}(0))\varsigma^{(n-i)}\left(0, \frac{1}{2}\right)$$

Using $\sum_{r=1}^{q} r\cos\left(\frac{2\pi rp}{q}\right)=\frac{q}{2}$, $p=1,2,...,q-1$ and letting $n=1$ in (3.23) gives us (3.14) again.

Letting $n=2$ in (3.21) gives us

$$\gamma_1\left(\frac{p}{q}\right)=\sum_{r=1}^{q}\cos\left(\frac{2\pi rp}{q}\right)\varsigma''\left(0,\frac{r}{q}\right)+\sum_{r=1}^{q}\left[2\mathbf{Y}_1(\mathbf{g}(0))\cos\left(\frac{2\pi rp}{q}\right)+\frac{\pi}{2}\sin\left(\frac{2\pi rp}{q}\right)\right]\varsigma'\left(0,\frac{r}{q}\right)$$

$$+\sum_{r=1}^{q}\left[\mathbf{Y}_2(\mathbf{g}(0))\cos\left(\frac{2\pi rp}{q}\right)+\pi\mathbf{Y}_1(\mathbf{h}(0))\sin\left(\frac{2\pi rp}{q}\right)\right]\varsigma\left(0,\frac{r}{q}\right)$$

$$=\sum_{r=1}^{q}\left[\varsigma''\left(0,\frac{r}{q}\right)-2[\gamma+\log(2\pi q)]\varsigma'\left(0,\frac{r}{q}\right)+\left[[\gamma+\log(2\pi q)]^2-\frac{1}{2}\varsigma(2)\right]\varsigma\left(0,\frac{r}{q}\right)\right]\cos\left(\frac{2\pi rp}{q}\right)$$

$$+\pi\sum_{r=1}^{q}\left[\varsigma'\left(0,\frac{r}{q}\right)+[\gamma+\log(2\pi q)]\varsigma\left(s,\frac{r}{q}\right)\right]\sin\left(\frac{2\pi rp}{q}\right)$$

and this may be written as

$$=\varsigma''(0)-\frac{1}{2}\left[[\gamma+\log(2\pi)]^2-\frac{1}{2}\varsigma(2)\right]-[\gamma+\log(2\pi)]\log q-\frac{1}{2}\log^2 q$$

$$+\sum_{r=1}^{q-1}\varsigma''\left(0,\frac{r}{q}\right)\cos\left(\frac{2\pi rp}{q}\right)-2[\gamma+\log(2\pi q)]\sum_{r=1}^{q-1}\log\Gamma\left(\frac{r}{q}\right)\cos\left(\frac{2\pi rp}{q}\right)$$

$$+\pi\sum_{r=1}^{q-1}\log\Gamma\left(\frac{r}{q}\right)\sin\left(\frac{2\pi rp}{q}\right)+\frac{\pi}{2}[\gamma+\log(2\pi q)]\cot\left(\frac{p\pi}{q}\right)$$

Letting $n=2$ $p=q=1$ in (3.21) gives us

$$\gamma_1=\varsigma''(0)+[\gamma+\log(2\pi)]\log(2\pi)-\frac{1}{2}\left([\gamma+\log(2\pi)]^2-\frac{1}{2}\varsigma(2)\right)$$

and hence we obtain

$$\gamma_1\left(\frac{p}{q}\right)=\gamma_1-[\gamma+\log(2\pi)]\log(2\pi q)-\frac{1}{2}\log^2 q$$

$$+\sum_{r=1}^{q-1}\varsigma''\left(0,\frac{r}{q}\right)\cos\left(\frac{2\pi rp}{q}\right)-2[\gamma+\log(2\pi q)]\sum_{r=1}^{q-1}\log\Gamma\left(\frac{r}{q}\right)\cos\left(\frac{2\pi rp}{q}\right)$$

$$+\pi\sum_{r=1}^{q-1}\log\Gamma\left(\frac{r}{q}\right)\sin\left(\frac{2\pi rp}{q}\right)+\frac{\pi}{2}[\gamma+\log(2\pi q)]\cot\left(\frac{p\pi}{q}\right)$$

Other equivalent representations for $\gamma_1\left(\frac{p}{q}\right)$ are contained in [7] and [10].

Some examples of the first Stieltjes constants are illustrated below:

$$\gamma_1\left(\frac{1}{4}\right)=\frac{1}{2}[2\gamma_1-7\log^2 2-6\gamma\log 2]-\frac{1}{2}\pi\left[\gamma+4\log 2+3\log\pi-4\log\Gamma\left(\frac{1}{4}\right)\right]$$

$$\gamma_1\left(\frac{3}{4}\right)=\frac{1}{2}[2\gamma_1-7\log^2 2-6\gamma\log 2]+\frac{1}{2}\pi\left[\gamma+4\log 2+3\log\pi-4\log\Gamma\left(\frac{1}{4}\right)\right]$$

$$\gamma_1\left(\frac{1}{5}\right)=\frac{1}{4}\left[4\gamma_1-\frac{5}{2}\log^2 5-5\gamma\log 5\right]-\frac{1}{2}\pi[\log(10\pi)+\gamma]\cot\left(\frac{\pi}{5}\right)$$

$$+\frac{1}{4}\sqrt{5}\left[\varsigma''\left(0,\frac{1}{5}\right)-\varsigma''\left(0,\frac{2}{5}\right)-\varsigma''\left(0,\frac{3}{5}\right)+\varsigma''\left(0,\frac{4}{5}\right)-[\gamma+\log(2\pi)]\log\frac{1}{2}(3+\sqrt{5})\right]$$

$$+\pi\left\{\sin\left(\frac{2\pi}{5}\right)\left[\log\Gamma\left(\frac{1}{5}\right)-\log\Gamma\left(\frac{4}{5}\right)\right]+\sin\left(\frac{\pi}{5}\right)\left[\log\Gamma\left(\frac{2}{5}\right)-\log\Gamma\left(\frac{3}{5}\right)\right]\right\}$$

**4. Some functional equations for $\varsigma''(0,x)$**

We note that the above representation of $\gamma_1\left(\frac{1}{5}\right)$ contains four terms involving the second derivative of the Hurwitz zeta function. As noted by Blagouchine [7], in certain cases it may be possible to reduce the number of terms by employing an appropriate functional equation for $\varsigma''(0,x)$.

The general Kubert identity is derived in [17, p.169]

(4.1) $$\Phi(s,x,z)=q^{-s}\sum_{r=0}^{q-1}\Phi\left(s,\frac{r+x}{q},z^q\right)$$

where $\Phi(s,x,z)$ is the Hurwitz-Lerch zeta function

(4.2) $$\Phi(s,x,z)=\sum_{n=0}^{\infty}\frac{z^n}{(n+x)^n}$$

We see that $\Phi(s,x,1)=\varsigma(s,x)$ and therefore we have

(4.3) $$q^s\varsigma(s,x)=\sum_{r=0}^{q-1}\varsigma\left(s,\frac{r+x}{q}\right)$$

Letting $x=1$ in (4.3) results in the well known identity

$$(q^s-1)\varsigma(s)=\sum_{r=1}^{q-1}\varsigma\left(s,\frac{r}{q}\right)$$

and differentiation gives us

$$\sum_{r=1}^{q-1}\varsigma''\left(0,\frac{r}{q}\right)=-\log q\log(2\pi)-\frac{1}{2}\log^2 q$$

Differentiation of (4.3) results in

(4.4) $$q^s\varsigma'(s,x)+q^s\varsigma(s,x)\log q=\sum_{r=0}^{q-1}\varsigma'\left(s,\frac{r+x}{q}\right)$$

and, letting $s=0$ and substituting Lerch's identity (3.12), we get

(4.5) $$\log\Gamma(x)=\sum_{r=0}^{q-1}\log\Gamma\left(\frac{r+x}{q}\right)-\frac{(q-1)}{2}\log(2\pi)-\left(\frac{1}{2}-x\right)\log q$$

or

(4.6) $$(2\pi)^{(q-1)/2}n^{(1/2)-x}\Gamma(x)=\prod_{r=0}^{q-1}\Gamma\left((r+x)/q\right)$$

which is Gauss's multiplication theorem for the gamma function [24, p.23]. The above derivation is contained in Milnor's paper [20].

Differentiation of (4.4) with respect to $s$ results in

(4.7) $$q^s\varsigma''(s,x)+2q^s\log q\varsigma'(s,x)+q^s\varsigma(s,x)\log^2 q=\sum_{r=0}^{q-1}\varsigma''\left(s,\frac{r+x}{q}\right)$$

and letting $s=0$ we obtain

(4.8) $$\varsigma''(0,x)+2\log q\varsigma'(0,x)+\varsigma(0,x)\log^2 q=\sum_{r=0}^{q-1}\varsigma''\left(0,\frac{r+x}{q}\right)$$

which may be expressed as [15]

(4.9) $$\varsigma''(0,x)+2\log q\left[\log\Gamma(x)-\frac{1}{2}\log(2\pi)\right]+\left[\frac{1}{2}-x\right]\log^2 q=\sum_{r=0}^{q-1}\varsigma''\left(0,\frac{r+x}{q}\right)$$

Letting $x \to 2x$ gives us

$$\varsigma''(0,2x) + 2\log q\left[\log\Gamma(2x) - \frac{1}{2}\log(2\pi)\right] + \left[\frac{1}{2} - 2x\right]\log^2 q = \sum_{r=0}^{q-1} \varsigma''\left(0, \frac{r+2x}{q}\right)$$

and with $q = 2$ this becomes

$$\varsigma''(0,2x) + 2\log 2\left[\log\Gamma(2x) - \frac{1}{2}\log(2\pi)\right] + \log^2 2\left[\frac{1}{2} - 2x\right] = \varsigma''(0,x) + \varsigma''\left(0, x+\frac{1}{2}\right)$$

or equivalently

$$\varsigma''(0,2x) - \varsigma''\left(0, x+\frac{1}{2}\right) = -2\log 2\left[\log\Gamma(2x) - \frac{1}{2}\log(2\pi)\right] - \log^2 2\left[\frac{1}{2} - 2x\right]$$

We also have the functional equation for 0 < $x$ < ½

(4.10) $$\varsigma''\left(0, x+\frac{1}{2}\right) + \varsigma''\left(0, \frac{1}{2} - x\right) = \varsigma''(0,2x) + \varsigma''(0,1-2x)$$

$$-\left[\varsigma''(0,x) + \varsigma''(0,1-x)\right] - 2\log 2\log(2\sin 2\pi x)$$

Differentiating (4.10) gives us the following functional equation involving the first Stieltjes constant:

(4.11) $$\gamma_1\left(x+\frac{1}{2}\right) - \gamma_1\left(\frac{1}{2} - x\right) = 2\left[\gamma_1(2x) - \gamma_1(1-2x)\right] - \left[\gamma_1(x) - \gamma_1(1-x)\right] - 2\pi\log 2\cdot\cot 2\pi x$$

We may also note that for $0 < x < 1$

(4.12) $$\sum_{n=1}^{\infty} \frac{\log n}{n}\cos 2\pi n x = \frac{1}{2}\left[\varsigma''(0,x) + \varsigma''(0,1-x)\right] + [\gamma + \log(2\pi)]\log(2\sin\pi x)$$

This identity was derived by Ramanujan and it appears, albeit in a highly disguised form, in Berndt's book [4, Part I, p.208] and [5]. This identity was reported by Deninger [15] in 1984 and was also derived by Blagouchine [7] in 2014.

## 5. A formula for the Stieltjes constants in terms of the higher derivatives of the Riemann zeta function $\varsigma^{(n)}(0)$.

We recall Hasse's formula (2.1) with $x = 1$

$$(s-1)\varsigma(s) = (s-1)\varsigma(s,1) = \sum_{n=0}^{\infty} \frac{1}{n+1} \sum_{k=0}^{n} \binom{n}{k} \frac{(-1)^k}{(1+k)^{s-1}}$$

and with $s \to 1-s$ we have

$$s\varsigma(1-s) = -\sum_{n=0}^{\infty}\frac{1}{n+1}\sum_{k=0}^{n}\binom{n}{k}(-1)^k(1+k)^s$$

We note the functional equation for the Riemann zeta function

$$2(2\pi)^{-s}\Gamma(s)\cos(\pi s/2)\varsigma(s) = \varsigma(1-s)$$

and, multiplying by $s$, we see that

$$f(s) = s\varsigma(1-s) = 2(2\pi)^{-s}\Gamma(s+1)\cos(\pi s/2)\varsigma(s) = -\sum_{n=0}^{\infty}\frac{1}{n+1}\sum_{k=0}^{n}\binom{n}{k}(-1)^k(1+k)^s$$

We have the $n$th derivative

$$f^{(n)}(s) = 2\sum_{i=0}^{n}\binom{n}{i}D^i[(2\pi)^{-s}\Gamma(s+1)]D^{n-i}[\cos(\pi s/2)\varsigma(s)]$$

and since

$$\frac{d}{ds}\Gamma(1+s)(2\pi)^{-s} = \Gamma(1+s)(2\pi)^{-s}[\psi(1+s)-\log(2\pi)]$$

we see from (3.6) that

$$D^i\Gamma(1+s)(2\pi)^{-s} = \Gamma(1+s)(2\pi)^{-s}Y_i\left(\psi(1+s)-\log(2\pi), \psi'(1+s), \cdots, \psi^{(i-1)}(1+s)\right)$$

and for convenience we designate $\mathbf{Y}_i(s)$ as

$$\mathbf{Y}_i(s) = Y_i\left(\psi(1+s)-\log(2\pi), \psi'(1+s), \cdots, \psi^{(i-1)}(1+s)\right)$$

We then have

$$D^i\Gamma(1+s)(2\pi)^{-s}\Big|_{s=0} = \mathbf{Y}_i(0)$$

where $\mathbf{Y}_i(0) = Y_i\left(-[\gamma+\log(2\pi)], \psi'(1), \cdots, \psi^{(i-1)}(1)\right)$.

Using (3.7) we then have

$$\mathbf{Y}_i(0) = Y_i\left(-[\gamma+\log(2\pi)], 1!\varsigma(2), \cdots, (-1)^i(i-1)!\varsigma(i)\right)$$

Accordingly, using (3.8), we may write

$$\mathbf{Y}_i(0) = (-1)^i Y_i\left(\gamma+\log(2\pi), 1!\varsigma(2), \cdots, (i-1)!\varsigma(i)\right)$$

which we abbreviate to $\mathbf{Y}_i(0) = (-1)^i \mathbf{Y}_i^*(0)$.

Reference to (2.4) shows that $f^{(n)}(0) = n\gamma_{n-1}$ and hence we have

(5.1) $$n\gamma_{n-1} = 2(-1)^n \sum_{i=0}^{n} \binom{n}{i} \mathbf{Y}_i^*(0) \sum_{j=0}^{n-i} \binom{n-i}{j} (-1)^j \left(\frac{\pi}{2}\right)^{n-i-j} \cos\left(\frac{\pi(n-i-j)}{2}\right) \varsigma^{(j)}(0)$$

This may be contrasted with the method employed by Apostol [2] who obtained the result:

$$\frac{(-1)^n \pi}{n!} \varsigma^{(n)}(0) = \frac{\operatorname{Im}\left(z^{n+1}\right)}{(n+1)!} + \sum_{k=1}^{n-1} a_k \frac{\operatorname{Im}\left(z^{n-k}\right)}{(n-k)!}$$

where $z = -\log(2\pi) - \frac{\pi}{2} i$ and the coefficients $a_k$ are given by the series

$$\Gamma(s)\varsigma(s) = \frac{1}{s-1} + \sum_{n=0}^{\infty} a_n (s-1)^n$$

An alternative, and more compact, representation of $\gamma_{n-1}$ may be derived as follows:

Using Euler's reflection formula

$$\Gamma(s)\Gamma(1-s) = \frac{\pi}{\sin \pi s}$$

we may express Riemann's functional equation as

(5.2) $$s\varsigma(1-s) = 2(2\pi)^{-s} \frac{\Gamma(1+s/2)\Gamma(1-s/2)}{\Gamma(1-s)} \varsigma(s)$$

For convenience we define $\lambda(s)$ by

$$\lambda(s) = (2\pi)^{-s} \frac{\Gamma(1+s/2)\Gamma(1-s/2)}{\Gamma(1-s)}$$

which differs slightly from the function used above in Section 3. We have

$$\lambda'(s) = \lambda(s) g(s)$$

where

$$g(s) = \psi(1-s) + \frac{1}{2}\psi\left(1+\frac{s}{2}\right) - \frac{1}{2}\psi\left(1-\frac{s}{2}\right) - \log(2\pi)$$

and hence

$$\lambda^{(i)}(s) = \lambda(s) Y_i\left(g(s), g^{(1)}(s), ..., g^{(i-1)}(s)\right)$$

where

$$g^{(r)}(s) = (-1)^r \psi^{(r)}(1-s) + \frac{1}{2^{r+1}}\psi^{(r)}\left(1+\frac{s}{2}\right) - \frac{1}{2^{r+1}}(-1)^r \psi^{(r)}\left(1-\frac{s}{2}\right)$$

and

$$g^{(r)}(0) = \left[\frac{1}{2^{r+1}}[(-1)^{r+1}+1]-1\right] r!\varsigma(r+1)$$

We see that

$$f^{(n)}(s) = 2\sum_{i=0}^{n}\binom{n}{i}\lambda^{(i)}(s)\varsigma^{(n-i)}(s)$$

and thus

$$f^{(n)}(0) = 2\sum_{i=0}^{n}\binom{n}{i} Y_i\left(g(0), g^{(1)}(0), ..., g^{(i-1)}(0)\right)\varsigma^{(n-i)}(0)$$

Since $f^{(n)}(0) = n\gamma_{n-1}$ we obtain

(5.3) $$n\gamma_{n-1} = 2\sum_{i=0}^{n}\binom{n}{i} Y_i\left(g(0), g^{(1)}(0), ..., g^{(i-1)}(0)\right)\varsigma^{(n-i)}(0)$$

which, for convenience, we express as

(5.4) $$n\gamma_{n-1} = 2\sum_{i=0}^{n}\binom{n}{i}\mathbf{Y}_i(\mathbf{g}(0))\varsigma^{(n-i)}(0)$$

where $\mathbf{Y}_i(\mathbf{g}(0)) \equiv Y_i\left(g(0), g^{(1)}(0), ..., g^{(i-1)}(0)\right)$.

For example, with $n = 2$ we have

$$\gamma_1 = \varsigma''(0) + 2\mathbf{Y}_1(\mathbf{g}(0))\varsigma'(0) + \mathbf{Y}_2(\mathbf{g}(0))\varsigma(0)$$

Noting that $\varsigma(0) = -\frac{1}{2}$ and, using Lerch's identity, we easily see that $\varsigma'(0) = -\frac{1}{2}\log(2\pi)$. We have

$$\mathbf{Y}_1(\mathbf{g}(0)) = -[\gamma + \log(2\pi)]$$

$$\mathbf{Y}_2(\mathbf{g}(0)) = [\gamma + \log(2\pi)]^2 - \frac{1}{2}\varsigma(2)$$

and thus we have

(5.5) $$\gamma_1 = \varsigma''(0) + [\gamma + \log(2\pi)]\log(2\pi) - \frac{1}{2}\left([\gamma + \log(2\pi)]^2 - \frac{1}{2}\varsigma(2)\right)$$

which corresponds with the result previously obtained by Ramanujan [5] and Apostol [2]

$$\varsigma''(0) = \gamma_1 + \frac{1}{2}\gamma^2 - \frac{1}{24}\pi^2 - \frac{1}{2}\log^2(2\pi)$$

**6. Miscellaneous results**

Li, Chakraborty and Kanemitsu [18] showed in 2010 that for integers $q > 1$ and $q > p$

(6.1)

$$q^{-s}\sum_{j=1}^{q}\exp\left(\frac{2\pi ijp}{q}\right)\varsigma\left(s,\frac{j}{q}\right) = i\frac{\Gamma(1-s)}{(2\pi)^{1-s}}\left[\exp\left(-\frac{\pi is}{2}\right)\varsigma\left(1-s,\frac{p}{q}\right) - \exp\left(\frac{\pi is}{2}\right)\varsigma\left(1-s,1-\frac{p}{q}\right)\right]$$

and, assuming that $s$ is real, we have the real and imaginary parts

(6.2) $$q^{-s}\sum_{j=1}^{q}\cos\left(\frac{2\pi jp}{q}\right)\varsigma\left(s,\frac{j}{q}\right) = \frac{\Gamma(1-s)}{(2\pi)^{1-s}}\sin\left(\frac{\pi s}{2}\right)\left[\varsigma\left(1-s,\frac{p}{q}\right) + \varsigma\left(1-s,1-\frac{p}{q}\right)\right]$$

(6.3) $$q^{-s}\sum_{j=1}^{q}\sin\left(\frac{2\pi jp}{q}\right)\varsigma\left(s,\frac{j}{q}\right) = \frac{\Gamma(1-s)}{(2\pi)^{1-s}}\cos\left(\frac{\pi s}{2}\right)\left[\varsigma\left(1-s,\frac{p}{q}\right) - \varsigma\left(1-s,1-\frac{p}{q}\right)\right]$$

We now employ these formulae to derive the following propositions which also appear in [7].

**Proposition 6.1**

(6.4) $$\sum_{j=1}^{q}\sin\left(\frac{2\pi jp}{q}\right)\psi\left(\frac{j}{q}\right) = \pi q\left[\frac{p}{q} - \frac{1}{2}\right]$$

**Proof**

Since $\sum_{j=1}^{q}\sin\left(\frac{2\pi jp}{q}\right) = 0$ we may write (6.3) as

(6.5) $$q^{-s}\sum_{j=1}^{q}\sin\left(\frac{2\pi jp}{q}\right)\left[\varsigma\left(s,\frac{j}{q}\right) - \varsigma(s)\right] = \frac{\Gamma(1-s)}{(2\pi)^{1-s}}\cos\left(\frac{\pi s}{2}\right)\left[\varsigma\left(1-s,\frac{p}{q}\right) - \varsigma\left(1-s,1-\frac{p}{q}\right)\right]$$

as noted by Chakraborty, Kanemitsu and Kuzumaki [8] in 2009.

Letting $s \to 1$ in (6.5) we obtain the limit

$$\lim_{s\to 1}\frac{\Gamma(1-s)}{(2\pi)^{1-s}}\cos\left(\frac{\pi s}{2}\right) = \lim_{s\to 1}(1-s)\Gamma(1-s)\frac{\cos\left(\frac{\pi s}{2}\right)}{1-s}$$

$$= \lim_{s\to 1}\Gamma(2-s)\lim_{s\to 1}\frac{\cos\left(\frac{\pi s}{2}\right)}{1-s}$$

and L'Hôpital's rule gives us

$$\lim_{s\to 1}\frac{\cos\left(\frac{\pi s}{2}\right)}{1-s}=\frac{\pi}{2}$$

Hence we obtain

$$q^{-1}\sum_{j=1}^{q}\sin\left(\frac{2\pi jp}{q}\right)\lim_{s\to 1}\left[\varsigma\left(s,\frac{j}{q}\right)-\varsigma(s)\right]=\frac{\pi}{2}\left[\varsigma\left(0,\frac{p}{q}\right)-\varsigma\left(0,1-\frac{p}{q}\right)\right]$$

so that

$$q^{-1}\sum_{j=1}^{q}\sin\left(\frac{2\pi jp}{q}\right)\left[\gamma_0\left(\frac{j}{q}\right)-\gamma\right]=\frac{\pi}{2}\left[1-\frac{2p}{q}\right]$$

Hence using (1.2) we have

$$\sum_{j=1}^{q}\sin\left(\frac{2\pi jp}{q}\right)\psi\left(\frac{j}{q}\right)=\pi q\left[\frac{p}{q}-\frac{1}{2}\right]$$

which is Gauss's second formula for the digamma function [24, p.19].

**Proposition 6.2**

(6.6) $$\sum_{j=1}^{q}\cos\left(\frac{2\pi jp}{q}\right)\psi\left(\frac{j}{q}\right)=q\log\left[2\sin\left(\frac{\pi p}{q}\right)\right]$$

**Proof**

We now consider (6.2) as $s\to 1$. We see that

$$\frac{\Gamma(1-s)}{(2\pi)^{1-s}}\sin\left(\frac{\pi s}{2}\right)\left[\varsigma\left(1-s,\frac{p}{q}\right)+\varsigma\left(1-s,1-\frac{p}{q}\right)\right]$$

$$=\frac{\pi}{(2\pi)^{1-s}\Gamma(s)\sin\pi s}\sin\left(\frac{\pi s}{2}\right)\left[\varsigma\left(1-s,\frac{p}{q}\right)+\varsigma\left(1-s,1-\frac{p}{q}\right)\right]$$

$$=\frac{\pi\left[\varsigma\left(1-s,\frac{p}{q}\right)+\varsigma\left(1-s,1-\frac{p}{q}\right)\right]}{(2\pi)^{1-s}\Gamma(s)2\cos\left(\frac{\pi s}{2}\right)}$$

Since $\varsigma(0,x)=\frac{1}{2}-x$ we see that $\varsigma(0,x)+\varsigma(0,1-x)=0$ and hence we may employ L'Hôpital's rule to show that

$$\lim_{s\to 1}\frac{\pi\left[\varsigma\left(1-s,\frac{p}{q}\right)+\varsigma\left(1-s,1-\frac{p}{q}\right)\right]}{(2\pi)^{1-s}\Gamma(s)2\cos\left(\frac{\pi s}{2}\right)}=\lim_{s\to 1}\frac{\varsigma'\left(1-s,\frac{p}{q}\right)+\varsigma'\left(1-s,1-\frac{p}{q}\right)}{\sin\left(\frac{\pi s}{2}\right)}$$

$$=\varsigma'\left(0,\frac{p}{q}\right)+\varsigma'\left(0,1-\frac{p}{q}\right)$$

Since $\varsigma'(0,x)+\varsigma'(0,1-x)=\log[\Gamma(x)\Gamma(1-x)]-\log(2\pi)$ the limit becomes $-\log\left[2\sin\left(\frac{\pi p}{q}\right)\right]$ and hence we obtain the second limb of Gauss's second formula for the digamma function [24, p.19]

$$\sum_{j=1}^{q}\cos\left(\frac{2\pi jp}{q}\right)\psi\left(\frac{j}{q}\right)=q\log\left[2\sin\left(\frac{\pi p}{q}\right)\right]$$

**Proposition 6.3**

(6.7)

$$\sum_{j=1}^{q}\sin\left(\frac{2\pi jp}{q}\right)\gamma_1\left(\frac{j}{q}\right)=\frac{\pi q}{2}\left[2\log\Gamma\left(\frac{p}{q}\right)-\log\pi+\log\sin\left(\frac{p\pi}{q}\right)\right]+\pi q\left[\gamma+\log(2\pi q)\right]\left[\frac{p}{q}-\frac{1}{2}\right]$$

**Proof**

Differentiating (6.5) gives us

(6.8)

$$\sum_{j=1}^{q}\sin\left(\frac{2\pi jp}{q}\right)\left[\varsigma'\left(s,\frac{j}{q}\right)-\varsigma'(s)\right]=-\frac{1}{2\pi}\Gamma(1-s)(2\pi q)^{s}\cos\left(\frac{\pi s}{2}\right)\left[\varsigma'\left(1-s,\frac{p}{q}\right)-\varsigma'\left(1-s,1-\frac{p}{q}\right)\right]$$

$$+\frac{1}{2\pi}\Gamma(1-s)(2\pi q)^{s}\cos\left(\frac{\pi s}{2}\right)\left[-\psi(1-s)+\log(2\pi q)-\frac{\pi}{2}\tan\left(\frac{\pi s}{2}\right)\right]\left[\varsigma\left(1-s,\frac{p}{q}\right)-\varsigma\left(1-s,1-\frac{p}{q}\right)\right]$$

With $s=1$ we obtain

$$\sum_{j=1}^{q}\sin\left(\frac{2\pi jp}{q}\right)\gamma_1\left(\frac{j}{q}\right)=\frac{\pi q}{2}\left[\varsigma'\left(0,\frac{p}{q}\right)-\varsigma'\left(0,1-\frac{p}{q}\right)\right]$$

$$-\frac{\pi q}{2}\left[\gamma+\log(2\pi q)\right]\left[\varsigma\left(0,\frac{p}{q}\right)-\varsigma\left(0,1-\frac{p}{q}\right)\right]$$

and using Lerch's identity this becomes

$$=\frac{\pi q}{2}\left[\log\Gamma\left(\frac{p}{q}\right)-\log\Gamma\left(1-\frac{p}{q}\right)\right]-\frac{\pi q}{2}[\gamma+\log(2\pi q)]\left[1-\frac{2p}{q}\right]$$

Substituting Euler's reflection formula

$$\log\Gamma\left(\frac{p}{q}\right)+\log\Gamma\left(1-\frac{p}{q}\right)=\log\pi-\log\sin\left(\frac{p\pi}{q}\right)$$

we then obtain the desired result

$$\sum_{j=1}^{q}\sin\left(\frac{2\pi jp}{q}\right)\gamma_1\left(\frac{j}{q}\right)=\frac{\pi q}{2}\left[2\log\Gamma\left(\frac{p}{q}\right)-\log\pi+\log\sin\left(\frac{p\pi}{q}\right)\right]-\frac{\pi q}{2}[\gamma+\log(2\pi q)]\left[1-\frac{2p}{q}\right]$$

**Proposition 6.4**

(6.9)

$$\sum_{j=1}^{q}\cos\left(\frac{2\pi jp}{q}\right)\gamma_1\left(\frac{j}{q}\right)=\frac{q}{2}\left[\varsigma''\left(0,\frac{p}{q}\right)+\varsigma''\left(0,1-\frac{p}{q}\right)\right]+q[\gamma+\log(2\pi)]\log\left[2\sin\left(\frac{p\pi}{q}\right)\right]$$

**Proof**

Since $\sum_{j=1}^{q}\cos\left(\frac{2\pi jp}{q}\right)=0$ we may write (6.2) as

$$\sum_{j=1}^{q}\cos\left(\frac{2\pi jp}{q}\right)\left[\varsigma\left(s,\frac{j}{q}\right)-\varsigma(s)\right]=\frac{\Gamma(1-s)}{2\pi}(2\pi q)^s\sin\left(\frac{\pi s}{2}\right)\left[\varsigma\left(1-s,\frac{p}{q}\right)+\varsigma\left(1-s,1-\frac{p}{q}\right)\right]$$

and differentiating this gives us

(6.10)

$$\sum_{j=1}^{q}\cos\left(\frac{2\pi jp}{q}\right)\left[\varsigma'\left(s,\frac{j}{q}\right)-\varsigma'(s)\right]=-\frac{1}{2\pi}\Gamma(1-s)(2\pi q)^s\sin\left(\frac{\pi s}{2}\right)\left[\varsigma'\left(1-s,\frac{p}{q}\right)+\varsigma'\left(1-s,1-\frac{p}{q}\right)\right]$$

$$+\frac{1}{2\pi}\Gamma(1-s)(2\pi q)^s\sin\left(\frac{\pi s}{2}\right)\left[-\psi(1-s)+\log(2\pi q)+\frac{\pi}{2}\cot\left(\frac{\pi s}{2}\right)\right]\left[\varsigma\left(1-s,\frac{p}{q}\right)+\varsigma\left(1-s,1-\frac{p}{q}\right)\right]$$

With $s=1$ we obtain

$$\sum_{j=1}^{q}\cos\left(\frac{2\pi jp}{q}\right)\gamma_1\left(\frac{j}{q}\right)$$

$$=-\lim_{s\to1}q\Gamma(1-s)\left\{\left[\varsigma'\left(0,\frac{p}{q}\right)+\varsigma'\left(0,1-\frac{p}{q}\right)\right]+[\psi(1-s)-\log(2\pi q)]\left[\varsigma\left(1-s,\frac{p}{q}\right)+\varsigma\left(1-s,1-\frac{p}{q}\right)\right]\right\}$$

and we will see below in Lemma 6.4 that defining $S$ by

$$S=\lim_{s\to 0}\Gamma(s)\left\{\left[\varsigma'(s,t)+\varsigma'(s,1-t)\right]+\left[\psi(s)-\log(2\pi)\right]\left[\varsigma(s,t)+\varsigma(s,1-t)\right]\right\}$$

we have

$$S=\frac{1}{2}\left[\varsigma''(0,t)+\varsigma''(0,1-t)\right]+[\gamma+\log(2\pi)]\log(2\sin\pi t)$$

Hence with $s\to 1-s$ we obtain

$$\sum_{j=1}^{q}\cos\left(\frac{2\pi jp}{q}\right)\gamma_1\left(\frac{j}{q}\right)=\frac{q}{2}\left[\varsigma''\left(0,\frac{p}{q}\right)+\varsigma''\left(0,1-\frac{p}{q}\right)\right]+q[\gamma+\log(2\pi)]\log\left[2\sin\left(\frac{p\pi}{q}\right)\right]$$

**Lemma 6.4**

$$\lim_{s\to 0}\Gamma(s)\left\{\left[\varsigma'(s,t)+\varsigma'(s,1-t)\right]+\left[\psi(s)-\log(2\pi)\right]\left[\varsigma(s,t)+\varsigma(s,1-t)\right]\right\}$$

$$=\frac{1}{2}\left[\varsigma''(0,t)+\varsigma''(0,1-t)\right]+[\gamma+\log(2\pi)]\log(2\sin\pi t)$$

**Proof**

As above, we designate $S$ as

$$S=\lim_{s\to 0}\Gamma(s)\left\{\left[\varsigma'(s,t)+\varsigma'(s,1-t)\right]+\left[\psi(s)-\log(2\pi)\right]\left[\varsigma(s,t)+\varsigma(s,1-t)\right]\right\}$$

$$=\lim_{s\to 0}s\Gamma(s)\left\{\frac{\left[\varsigma'(s,t)+\varsigma'(s,1-t)\right]+\left[\psi(s)-\log(2\pi)\right]\left[\varsigma(s,t)+\varsigma(s,1-t)\right]}{s}\right\}$$

(6.11) $$=\lim_{s\to 0}\left\{\frac{\left[\varsigma'(s,t)+\varsigma'(s,1-t)\right]+\left[\psi(s)-\log(2\pi)\right]\left[\varsigma(s,t)+\varsigma(s,1-t)\right]}{s}\right\}$$

We may write part of the numerator of (6.11) as

$$\psi(s)\left[\varsigma(s,t)+\varsigma(s,1-t)\right]=\frac{s\psi(s)\left[\varsigma(s,t)+\varsigma(s,1-t)\right]}{s}$$

and thus

$$\lim_{s\to 0}\psi(s)\left[\varsigma(s,t)+\varsigma(s,1-t)\right]=\lim_{s\to 0}s\psi(s)\lim_{s\to 0}\frac{\left[\varsigma(s,t)+\varsigma(s,1-t)\right]}{s}$$

Since $\psi(1+s)=\psi(s)+1/s$, we see that $\lim_{s\to 0}s\psi(s)=-1$ and, using L'Hôpital's rule, we obtain

$$\lim_{s\to 0}\psi(s)\left[\varsigma(s,t)+\varsigma(s,1-t)\right]=-\left[\varsigma'(0,t)+\varsigma'(0,1-t)\right]$$

Therefore we see that

$$\lim_{s\to 0}\left[\varsigma'(s,t)+\varsigma'(s,1-t)\right]+\left[\psi(s)-\log(2\pi)\right]\left[\varsigma(s,t)+\varsigma(s,1-t)\right]=0$$

and hence using L'Hôpital's rule at (6.11) we have

$$S=\varsigma''(0,t)+\varsigma''(0,1-t)+\lim_{s\to 0}\frac{d}{ds}\left[\psi(s)-\log(2\pi)\right]\left[\varsigma(s,t)+\varsigma(s,1-t)\right]$$

This limit is

$$=\lim_{s\to 0}\left[\psi(s)-\log(2\pi)\right]\left[\varsigma'(s,t)+\varsigma'(s,1-t)\right]+\psi'(s)\left[\varsigma(s,t)+\varsigma(s,1-t)\right]$$

$$=\lim_{s\to 0}\psi(s)\left[\varsigma'(s,t)+\varsigma'(s,1-t)\right]+\psi'(s)\left[\varsigma(s,t)+\varsigma(s,1-t)\right]+\log(2\pi)\log(2\sin\pi t)$$

$$=\lim_{s\to 0}\left[\psi(1+s)-\frac{1}{s}\right]\left[\varsigma'(s,t)+\varsigma'(s,1-t)\right]+\left[\psi'(1+s)+\frac{1}{s^2}\right]\left[\varsigma(s,t)+\varsigma(s,1-t)\right]$$

$$+\log(2\pi)\log(2\sin\pi t)$$

$$=-\gamma\left[\varsigma'(0,t)+\varsigma'(0,1-t)\right]-\lim_{s\to 0}\frac{1}{s}\left[\varsigma'(s,t)+\varsigma'(s,1-t)\right]+\frac{1}{s^2}\left[\varsigma(s,t)+\varsigma(s,1-t)\right]$$

$$+\log(2\pi)\log(2\sin\pi t)$$

$$=-\lim_{s\to 0}\frac{1}{s}\left[\varsigma'(s,t)+\varsigma'(s,1-t)\right]+\frac{1}{s^2}\left[\varsigma(s,t)+\varsigma(s,1-t)\right]+[\gamma+\log(2\pi)]\log(2\sin\pi t)$$

$$=\lim_{s\to 0}\frac{\left[\varsigma(s,t)+\varsigma(s,1-t)\right]-s\left[\varsigma'(s,t)+\varsigma'(s,1-t)\right]}{s^2}+[\gamma+\log(2\pi)]\log(2\sin\pi t)$$

Applying L'Hôpital's rule again gives us

$$\lim_{s\to 0}\frac{\left[\varsigma(s,t)+\varsigma(s,1-t)\right]-s\left[\varsigma'(s,t)+\varsigma'(s,1-t)\right]}{s^2}=-\lim_{s\to 0}\frac{s\left[\varsigma''(s,t)+\varsigma''(s,1-t)\right]}{2s}$$

$$=-\frac{1}{2}\left[\varsigma''(0,t)+\varsigma''(0,1-t)\right]$$

The lemma therefore follows.

**Proposition 6.5**

(6.12) $$\sum_{j=1}^{q-1}\sin\left(\frac{2\pi jp}{q}\right)\log\Gamma\left(\frac{j}{q}\right)=\frac{1}{2\pi}\left[\gamma_1\left(\frac{p}{q}\right)-\gamma_1\left(1-\frac{p}{q}\right)\right]+\frac{1}{2}[\gamma+\log(2\pi q)]\cot\left(\frac{p\pi}{q}\right)$$

**Proof**

With $s=0$ in (6.8) we obtain

$$\sum_{j=1}^{q}\sin\left(\frac{2\pi jp}{q}\right)\varsigma'\left(0,\frac{j}{q}\right)=-\frac{1}{2\pi}\left[\varsigma'\left(1,\frac{p}{q}\right)-\varsigma'\left(1,1-\frac{p}{q}\right)\right]$$

$$+\frac{1}{2\pi}[\gamma+\log(2\pi q)]\left[\varsigma\left(1,\frac{p}{q}\right)-\varsigma\left(1,1-\frac{p}{q}\right)\right]$$

$$=\frac{1}{2\pi}\left[\gamma_1\left(\frac{p}{q}\right)-\gamma_1\left(1-\frac{p}{q}\right)\right]$$

$$+\frac{1}{2\pi}[\gamma+\log(2\pi q)]\left[\gamma_0\left(\frac{p}{q}\right)-\gamma_0\left(1-\frac{p}{q}\right)\right]$$

$$=\frac{1}{2\pi}\left[\gamma_1\left(\frac{p}{q}\right)-\gamma_1\left(1-\frac{p}{q}\right)\right]$$

$$-\frac{1}{2\pi}[\gamma+\log(2\pi q)]\left[\psi\left(\frac{p}{q}\right)-\psi\left(1-\frac{p}{q}\right)\right]$$

and hence we obtain

(6.13) $$\sum_{j=1}^{q}\sin\left(\frac{2\pi jp}{q}\right)\varsigma'\left(0,\frac{j}{q}\right)=\frac{1}{2\pi}\left[\gamma_1\left(\frac{p}{q}\right)-\gamma_1\left(1-\frac{p}{q}\right)\right]+\frac{1}{2}[\gamma+\log(2\pi q)]\cot\left(\frac{p\pi}{q}\right)$$

We may also write this as

$$\sum_{j=1}^{q-1}\sin\left(\frac{2\pi jp}{q}\right)\log\Gamma\left(\frac{j}{q}\right)=\frac{1}{2\pi}\left[\gamma_1\left(\frac{p}{q}\right)-\gamma_1\left(1-\frac{p}{q}\right)\right]+\frac{1}{2}[\gamma+\log(2\pi q)]\cot\left(\frac{p\pi}{q}\right)$$

as originally derived by Malmstén [7] in 1849.

**Proposition 6.6**

(6.14) $$\sum_{j=1}^{q}\cos\left(\frac{2\pi jp}{q}\right)\log\Gamma\left(\frac{j}{q}\right)=-\frac{1}{4}\left[\psi\left(\frac{p}{q}\right)+\psi\left(1-\frac{p}{q}\right)\right]-\frac{1}{2}[\gamma+\log(2\pi q)]$$

**Proof**

With $s=0$ in (6.10) we obtain

$$2\pi\sum_{j=1}^{q}\cos\left(\frac{2\pi jp}{q}\right)\varsigma'\left(0,\frac{j}{q}\right)$$

$$=\lim_{s\to 0}\left(\frac{\pi}{2}\cos\left(\frac{\pi s}{2}\right)\left[\varsigma\left(1-s,\frac{p}{q}\right)+\varsigma\left(1-s,1-\frac{p}{q}\right)\right]-\sin\left(\frac{\pi s}{2}\right)\left[\varsigma'\left(1-s,\frac{p}{q}\right)+\varsigma'\left(1-s,1-\frac{p}{q}\right)\right]\right)$$

$$+\lim_{s\to 0}\sin\left(\frac{\pi s}{2}\right)\left[-\psi(1-s)+\log(2\pi q)\right]\left[\varsigma\left(1-s,\frac{p}{q}\right)+\varsigma\left(1-s,1-\frac{p}{q}\right)\right]$$

We express the first limit as

$$\lim_{s\to 0}\left(\frac{\pi}{2}\cos\left(\frac{\pi s}{2}\right)\left[\varsigma\left(1-s,\frac{p}{q}\right)+\varsigma\left(1-s,1-\frac{p}{q}\right)\right]-\sin\left(\frac{\pi s}{2}\right)\left[\varsigma'\left(1-s,\frac{p}{q}\right)+\varsigma'\left(1-s,1-\frac{p}{q}\right)\right]\right)$$

$$=\lim_{s\to 0}\frac{\pi}{2}\cos\left(\frac{\pi s}{2}\right)\left[\varsigma\left(1-s,\frac{p}{q}\right)+\varsigma\left(1-s,1-\frac{p}{q}\right)+\frac{2}{s}\right]$$

$$-\lim_{s\to 0}\sin\left(\frac{\pi s}{2}\right)\left[\varsigma'\left(1-s,\frac{p}{q}\right)+\varsigma'\left(1-s,1-\frac{p}{q}\right)+\frac{2}{s^2}\right]$$

$$-\lim_{s\to 0}\left[\frac{\pi}{2}\frac{2}{s}\cos\left(\frac{\pi s}{2}\right)-\frac{2}{s^2}\sin\left(\frac{\pi s}{2}\right)\right]$$

We note from (1.1) that

$$\varsigma(1-s,u)+\frac{1}{s}=\sum_{p=0}^{\infty}\frac{\gamma_p(u)}{p!}s^p$$

so that

$$\varsigma'(1-s,u)+\frac{1}{s^2}=-\sum_{p=0}^{\infty}\frac{\gamma_p(u)}{p!}ps^{p-1}$$

and we see that

$$\lim_{s\to 0}\frac{\pi}{2}\cos\left(\frac{\pi s}{2}\right)\left[\varsigma\left(1-s,\frac{p}{q}\right)+\varsigma\left(1-s,1-\frac{p}{q}\right)+\frac{2}{s}\right]=\frac{\pi}{2}\left[\gamma_0\left(\frac{p}{q}\right)+\gamma_0\left(1-\frac{p}{q}\right)\right]$$

Similarly we see that

$$\lim_{s\to 0}\sin\left(\frac{\pi s}{2}\right)\left[\varsigma'\left(1-s,\frac{p}{q}\right)+\varsigma'\left(1-s,1-\frac{p}{q}\right)+\frac{2}{s^2}\right]=0$$

We write

$$\lim_{s\to 0}\left[\frac{\pi}{2}\frac{2}{s}\cos\left(\tfrac{\pi s}{2}\right)-\frac{2}{s^2}\sin\left(\tfrac{\pi s}{2}\right)\right]=\lim_{s\to 0}\left[\frac{\pi s\cos\left(\tfrac{\pi s}{2}\right)-2\sin\left(\tfrac{\pi s}{2}\right)}{s^2}\right]$$

and applying L'Hôpital's rule we see that this limit vanishes.

We have the second limit

$$\lim_{s\to 0}\sin\left(\frac{\pi s}{2}\right)\left[-\psi(1-s)+\log(2\pi q)\right]\left[\varsigma\left(1-s,\frac{p}{q}\right)+\varsigma\left(1-s,1-\frac{p}{q}\right)\right]$$

$$=\lim_{s\to 0}\frac{\pi}{2}\frac{\sin\left(\tfrac{\pi s}{2}\right)}{\tfrac{\pi s}{2}}\left[-\psi(1-s)+\log(2\pi q)\right]s\left[\varsigma\left(1-s,\frac{p}{q}\right)+\varsigma\left(1-s,1-\frac{p}{q}\right)\right]$$

$$=-\pi\left[\gamma+\log(2\pi q)\right]$$

since $\lim_{s\to 0}s\varsigma(1-s,x)=-1$.
Hence we have obtained

$$\sum_{j=1}^{q}\cos\left(\frac{2\pi jp}{q}\right)\varsigma'\left(0,\frac{j}{q}\right)=\frac{1}{4}\left[\gamma_0\left(\frac{p}{q}\right)+\gamma_0\left(1-\frac{p}{q}\right)\right]-\frac{1}{2}[\gamma+\log(2\pi q)]$$

Substituting Lerch's identity we obtain a Gauss-type identity

$$\sum_{j=1}^{q}\cos\left(\frac{2\pi jp}{q}\right)\log\Gamma\left(\frac{j}{q}\right)=-\frac{1}{4}\left[\psi\left(\frac{p}{q}\right)+\psi\left(1-\frac{p}{q}\right)\right]-\frac{1}{2}[\gamma+\log(2\pi q)]$$

and we may write this in the form originally discovered by Malmstén (see [3], [19] and also Blagouchine's paper [6])

(6.15) $$2\sum_{j=1}^{q}\cos\left(\frac{2\pi jp}{q}\right)\log\Gamma\left(\frac{j}{q}\right)=-\frac{\pi}{2}\cot\left(\frac{\pi p}{q}\right)-\left[\gamma+\log(2\pi q)\right]-\psi\left(\frac{p}{q}\right)$$

This analysis may be continued to obtain formulae for $\sum_{j=1}^{q}\sin\left(\frac{2\pi jp}{q}\right)\varsigma^{(n)}\left(0,\frac{j}{q}\right)$ and $\sum_{j=1}^{q}\sin\left(\frac{2\pi jp}{q}\right)\varsigma^{(n)}\left(0,\frac{j}{q}\right)$; for example we have:

**Proposition 6.7**

(6.16) $$2\pi\sum_{j=1}^{q}\sin\left(\frac{2\pi jp}{q}\right)\varsigma''\left(0,\frac{j}{q}\right)=\gamma_2\left(\frac{p}{q}\right)-\gamma_2\left(1-\frac{p}{q}\right)$$

$$+4\pi\left[\gamma+\log(2\pi q)\right]\sum_{j=1}^{q}\sin\left(\frac{2\pi jp}{q}\right)\varsigma'\left(0,\frac{j}{q}\right)-\pi\left(\left[\gamma+\log(2\pi q)\right]^2+\frac{\pi^2}{12}\right)\cot\left(\frac{p\pi}{q}\right)$$

**Proof**

We have the derivative of (6.13)

$$\sum_{j=1}^{q}\sin\left(\frac{2\pi jp}{q}\right)\varsigma''\left(s,\frac{j}{q}\right)=\frac{1}{2\pi}\Gamma(1-s)(2\pi q)^s\cos\left(\frac{\pi s}{2}\right)\left[\varsigma''\left(1-s,\frac{p}{q}\right)-\varsigma''\left(1-s,1-\frac{p}{q}\right)\right]$$

$$-\frac{1}{\pi}\Gamma(1-s)(2\pi q)^s\cos\left(\frac{\pi s}{2}\right)\left[-\psi(1-s)+\log(2\pi q)-\frac{\pi}{2}\tan\left(\frac{\pi s}{2}\right)\right]\left[\varsigma'\left(1-s,\frac{p}{q}\right)-\varsigma'\left(1-s,1-\frac{p}{q}\right)\right]$$

$$+\frac{1}{2\pi}\Gamma(1-s)(2\pi q)^s\cos\left(\frac{\pi s}{2}\right)\left[-\psi(1-s)+\log(2\pi q)-\frac{\pi}{2}\tan\left(\frac{\pi s}{2}\right)\right]^2\left[\varsigma\left(1-s,\frac{p}{q}\right)-\varsigma\left(1-s,1-\frac{p}{q}\right)\right]$$

$$+\frac{1}{2\pi}\Gamma(1-s)(2\pi q)^s\cos\left(\frac{\pi s}{2}\right)\left[\psi'(1-s)+\log(2\pi q)-\frac{\pi^2}{4}\sec^2\left(\frac{\pi s}{2}\right)\right]\left[\varsigma\left(1-s,\frac{p}{q}\right)-\varsigma\left(1-s,1-\frac{p}{q}\right)\right]$$

With $s=0$ we obtain

$$\sum_{j=1}^{q}\sin\left(\frac{2\pi jp}{q}\right)\varsigma''\left(0,\frac{j}{q}\right)=\frac{1}{2\pi}\left[\gamma_2\left(\frac{p}{q}\right)-\gamma_2\left(1-\frac{p}{q}\right)\right]$$

$$+\frac{1}{\pi}\left[\gamma+\log(2\pi q)\right]\left[\gamma_1\left(\frac{p}{q}\right)-\gamma_1\left(1-\frac{p}{q}\right)\right]$$

$$+\frac{1}{2}\left[\gamma+\log(2\pi q)\right]^2\cot\left(\frac{p\pi}{q}\right)$$

$$+\frac{1}{2}\left[-\frac{\pi^2}{12}\right]\cot\left(\frac{p\pi}{q}\right)$$

Using (6.13) we obtain

$$2\pi\sum_{j=1}^{q}\sin\left(\frac{2\pi jp}{q}\right)\varsigma''\left(0,\frac{j}{q}\right)=\gamma_2\left(\frac{p}{q}\right)-\gamma_2\left(1-\frac{p}{q}\right)$$

$$+4\pi\left[\gamma+\log(2\pi q)\right]\sum_{j=1}^{q}\sin\left(\frac{2\pi jp}{q}\right)\varsigma'\left(0,\frac{j}{q}\right)-\pi\left(\left[\gamma+\log(2\pi q)\right]^2+\frac{\pi^2}{12}\right)\cot\left(\frac{p\pi}{q}\right)$$

or alternatively

$$\gamma_2\left(\frac{p}{q}\right)-\gamma_2\left(1-\frac{p}{q}\right)=2\pi\sum_{j=1}^{q-1}\sin\left(\frac{2\pi jp}{q}\right)\varsigma''\left(0,\frac{j}{q}\right)-4\pi[\gamma+\log(2\pi q)]\sum_{j=1}^{q}\sin\left(\frac{2\pi jp}{q}\right)\log\Gamma\left(\frac{j}{q}\right)$$

$$+\frac{\pi}{2}\left(2[\gamma+\log(2\pi q)]^2+\frac{\pi^2}{6}\right)\cot\left(\frac{\pi p}{q}\right)$$

**7. Open access to our own work**

This paper contains references to a number of other papers and, fortunately, most of them are currently freely available on the internet. Surely now is the time that all of our work should be freely accessible by all. The mathematics community should lead the way on this by publishing everything on arXiv, or in an equivalent open access repository. We think it, we write it, so why hide it? You know it makes sense.

Wessex House,
Devizes Road,
Upavon,
Pewsey,
Wiltshire SN9 6DL
dconnon@btopenworld.com